\documentclass[12pt,centertags,oneside]{amsart}

\usepackage{amsmath,amstext,amsthm,amscd,typearea,hyperref,stmaryrd,mathabx}
\usepackage{amssymb}
\usepackage{a4wide}
\usepackage[mathscr]{eucal}
\usepackage{mathrsfs}
\usepackage{typearea}
\usepackage{charter}
\usepackage{pdfsync}
\usepackage[a4paper,width=16.2cm,top=3cm,bottom=3cm,lines=50]{geometry}
\usepackage[all]{xy}
\usepackage{tikz}
\usepackage{enumitem}
\usepackage[foot]{amsaddr}
\usepackage{academicons}
\usepackage{xcolor}

\numberwithin{equation}{section}
\newcommand{\orcid}[1]{\href{https://orcid.org/#1}{\textsc{orc}i\textsc{d}}}

\title[]{Heat diffusions on holomorphic foliations with non-degenerate singularities}

\author{Fran\c cois Bacher}
\address{Universit\'e de Lille, 
Laboratoire de math\'ematiques Paul Painlev\'e, 
CNRS U.M.R. 8524,  
59655 Villeneuve d'Ascq Cedex, 
France.}
\email{francois.bacher@univ-lille.fr}
\thanks{\textsc{orc}i\textsc{d}: 0000-0001-8212-7109}

\date{April 1, 2023.}

\def\restriction#1#2{\mathchoice
              {\setbox1\hbox{${\displaystyle #1}_{\scriptstyle #2}$}
              \restrictionaux{#1}{#2}}
              {\setbox1\hbox{${\textstyle #1}_{\scriptstyle #2}$}
              \restrictionaux{#1}{#2}}
              {\setbox1\hbox{${\scriptstyle #1}_{\scriptscriptstyle #2}$}
              \restrictionaux{#1}{#2}}
              {\setbox1\hbox{${\scriptscriptstyle #1}_{\scriptscriptstyle #2}$}
              \restrictionaux{#1}{#2}}}
\def\restrictionaux#1#2{{#1\,\smash{\vrule height .8\ht1 depth .85\dp1}}_{\,#2}}

\theoremstyle{plain}
\newtheorem{thm}{Theorem}[section]
\newtheorem{lem}[thm]{Lemma}

\newtheorem{prop}[thm]{Proposition}

\newtheorem*{thm*}{Theorem}
\newtheorem*{conj*}{Conjecture}

\theoremstyle{definition}
\newtheorem{defn}[thm]{Definition}

\newtheorem*{exmp*}{Example}

\theoremstyle{remark}
\newtheorem{rem}[thm]{Remark}

\newtheorem*{pr}{Proof}

\makeatletter

\@addtoreset{equation}{section}
\makeatother

\DeclareMathOperator{\id}{id}

\DeclareMathOperator{\Dom}{Dom}

\DeclareMathOperator{\supp}{supp}

\DeclareMathOperator{\Leb}{Leb}

\newcommand{\cjg}[1]{\overline{#1}}

\newcommand{\adh}[1]{\overline{#1}}

\newcommand{\PC}{P}

\newcommand{\eps}{\varepsilon}

\newcommand{\fol}{\mathscr{F}}

\newcommand{\leafatlas}{\mathscr{L}}

\newcommand{\plfol}{\left(\mani{M},\leafatlas,\mani{E}\right)}

\newcommand{\zlogz}[1]{#1\log^{\star}#1}

\newcommand{\Lap}{\Delta}

\newcommand{\LapP}{\Lap_{\PC}}

\newcommand{\LapPv}[1]{\Lap_{\PC,#1}}

\newcommand{\tLap}{\tilde{\Lap}}

\newcommand{\tLapP}{\tLap_{\PC}}

\newcommand{\gr}{\nabla}

\newcommand{\ddc}{dd^{\text{c}}}

\newcommand{\hd}{\partial}

\newcommand{\bhd}{\cjg{\partial}}

\newcommand{\set}[1]{\mathbb{#1}}

\newcommand{\der}[2]{\frac{\partial#1}{\partial#2}}

\newcommand{\leaf}{L}

\newcommand{\leafu}[1]{\leaf_{#1}}

\newcommand{\norm}[1]{\left\Vert#1\right\Vert}

\newcommand{\intcc}[2]{\left[#1,#2\right]}

\newcommand{\intoo}[2]{\left(#1,#2\right)}

\newcommand{\intoc}[2]{\left(#1,#2\right]}

\newcommand{\Cmod}[1]{\left\vert#1\right\vert}

\newcommand{\indic}[1]{\mathbf{1}_{#1}}

\newcommand{\proj}[1]{\set{P}^{#1}}

\newcommand{\class}[1]{\mathscr{C}^{#1}}

\newcommand{\mani}[1]{#1}

\newcommand{\manis}[2]{\mani{#1}\backslash\mani{#2}}

\newcommand{\dsing}[3]{d_{\mani{#1}}(#3,\mani{#2})}

\newcommand{\rD}[1]{#1\set{D}}

\newcommand{\DR}[1]{\set{D}_{#1}}

\newcommand{\setst}{~\vert~}

\newcommand{\wo}[2]{#1\backslash#2}

\newcommand{\testfun}[1]{\mathscr{D}\left(#1\right)}

\newcommand{\testfunfol}{\testfun{\fol}}

\newcommand{\testformfol}[1]{\mathscr{D}^{#1}\left(\fol{}\right)}

\newcommand{\dherm}[3]{d_{\mani{#1}}(#2,#3)}

\newcommand{\dPC}[2]{d_{\PC}(#1,#2)}

\newcommand{\upb}{F}

\newcommand{\metPC}{g_{\PC}}

\newcommand{\metm}[1]{g_{\mani{#1}}}

\newcommand{\metmp}[2]{g_{\mani{#1},#2}}

\newcommand{\textlam}{Riemann surface lamination}

\newcommand{\textfol}{holomorphic foliation}

\newcommand{\textsingfol}{singular \textfol{}}

\newcommand{\textsinglams}{\textlam{}s with singularities}

\newcommand{\textsingfols}{\textsingfol{}s}

\newcommand{\texthypsingfol}{hyperbolic \textsingfol{}}

\begin{document}

\theoremstyle{plain}

\begin{abstract} Consider a Brody hyperbolic foliation with non-degenerate singularities on a compact complex manifold. We show that the leafwise heat diffusions and the abstract heat diffusions coincide. In particular, this will imply that the abstract heat diffusions are unique.
  \end{abstract}

\maketitle

\section{Introduction}

The global dynamics of \textlam{}s have recently received much attention. Much progress has been focused on developing an adapted ergodic theory. This could be a powerful tool to understand the global behaviour of \textlam{}s. The case of the projective space is maybe the most typical, since polynomial vector fields can be compactified naturally into foliations on $\proj{n}$. Lins Neto and Soares in~\cite{LNS} and Jouanolou in~\cite{Jou} have shown that a generic foliation $\fol{}$ of a given degree $d\geq2$ on $\proj{n}$ has only non-degenerate singularities. By a result of Lins Neto~\cite{LN2} and Glutsyuk~\cite{Glu}, a foliation having only non-degenerate singularities is necessarily hyperbolic, and even Brody hyperbolic in the sense of~\cite{DNSII}. Therefore, many authors have studied hyperbolic \textsinglams{} and in particular hyperbolic \textsingfols{}. We present briefly the results of some recent work and refer the reader to the surveys~\cite{surDinhSib,surForSib,surVANG18,surVANG21} for more details.

To set up an ergodic theory, one needs at least some notions of time, time average, space average (\emph{i.e.} measure), invariant and ergodic measure. It turns out that for foliations, it will also be convenient to use currents instead of measures and forms instead of functions. In the case of hyperbolic foliations, the Poincar\'{e} distance in a universal covering can be seen as a canonical time. As an analogue to invariant measures, it is natural to consider directed closed currents. This has led to various results and dynamical methods (see for example Rebelo~\cite{Reb}), but a large class of foliations do not have any. That is why Garnett has considered the weaker notions of harmonic currents and harmonic measures in~\cite{Gar}. This approach has led to strong results, see for example the survey~\cite{surForSib} of Forn\ae{}ss and Sibony or the unique ergodicity result of Dinh, Nguy\^{e}n, Sibony in~\cite{DNSuniergo}. The notion of quasi-harmonic measure was also introduced by Nguy\^{e}n in~\cite{NguOse} to establish an Oseledec multiplicative ergodic theorem for laminations. Then, the natural processes to take the average of a function on the leaves (that is, on the orbits) becomes the heat diffusions. Two approaches to define such diffusions have been tried and both have given a series of ergodic theorems. It is natural to wonder whether these two approaches lead to the same object, that is, whether both heat diffusions coincide. Though the equations they satisfy look similar, it is not clear and not proven if this is true in general. In his survey~\cite{surVANG21}, Nguy\^{e}n has worked on this question. His result is the following.

\begin{thm}[{\cite[Corollary~5.23]{surVANG21}}] \label{impcorVANG} Let $\fol{}=\plfol{}$ be a Brody hyperbolic compact \textsingfol{}. Suppose that all the singularities of $\fol{}$ are linearizable and hyperbolic. Then, for every harmonic measure that does not give mass to any leaf, the abstract heat diffusions and the leafwise heat diffusions coincide. In particular, the abstract heat diffusions are unique, \emph{i.e.}, they do not depend on the considered harmonic measures.
\end{thm}

He mentions as a question (see~\cite[Problem~5.24]{surVANG21}) to find better sufficient conditions for the two heat diffusions to coincide. Our main result gives an improvement of Theorem~\ref{impcorVANG}.

\begin{thm}\label{impcor} Let $\fol=\plfol{}$ be a Brody hyperbolic compact \textsingfol{}. Suppose that all the singularities of $\fol{}$ are non-degenerate. Then, for every harmonic measure, the abstract heat diffusions coincide with the leafwise heat diffusions. In particular, the abstract heat diffusions are unique, \emph{i.e.}, they are independent of the harmonic measures. 
\end{thm}

Let us explain the method of our proof. When we study non-degenerate singularities, we encounter more singular objects than in the case of linearizable weakly hyperbolic singularities. In our context, the current associated to the harmonic measure may have non-vanishing Lelong numbers on the singularities. Indeed, Chen has shown in~\cite{Chen} that the Lelong number of a directed harmonic current on a non-hyperbolic singularity may be positive. On the other hand, Nguy\^{e}n has shown in~\cite{VANGLelong} that it cannot if the singularity is linearizable and weakly hyperbolic. To deal with this new difficulty, we have to control the heat diffusions near the singularities. More precisely, given an initial heat distribution far from the singular set, we show that it essentially does not reach the neighbourhoods of a non-degenerate singularity in finite time (see Lemma~\ref{factepsinf}). This result relies on estimates of the Poincar\'{e} metric near the singularities.

The article is organized as follows. In Section~\ref{secPC}, we introduce the leafwise Poincar\'{e} metric $\metPC{}$. We specify the type of singularities we will consider. In the case where a foliation $\fol{}=\plfol{}$ is also endowed with a Riemannian metric $\metm{M}$, we define the function $\eta$. It is a quantitative way to compare $\metPC{}$ and $\metm{M}$ so that $\eta^2\metPC{}=4\metm{M}$. This is the main tool we use to estimate the diffusions near the singularities. In Section~\ref{secDt}, we describe the heat kernel and the leafwise heat diffusions, as well as their behaviours with respect to a uniformization of a leaf. In Section~\ref{seccur}, we explain the link between harmonic measures and harmonic currents. In our context, we show that the Poincar\'{e} mass of a directed harmonic current is finite. In Section~\ref{secSt}, we recall the construction of the abstract heat diffusions. We state some intermediate lemmas of~\cite{DNS12} to show its existence, because we will need them for our proof. Section~\ref{secidheats} is devoted to the proof of Theorem~\ref{impcor}. We state and prove an abstract criterion, analogous to~\cite[Theorem~5.17]{surVANG21} in Theorem~\ref{mainthm}. In most of the article, we follow the method and notations of~\cite{surVANG21}.

\subsection*{Notations} Throughout this paper, we will denote by $\set{D}$ the unit disk of $\set{C}$ and $\set{B}$ the unit ball in $\set{C}^k$. We denote by $\rD{r}$ the disk of radius $r$ in $\set{C}$ and $r\set{B}$ the ball of radius $r$ in $\set{C}^k$. When we consider the hyperbolic distance in $\set{D}$, we denote by $\DR{R}$ the disk of hyperbolic radius $R$ in $\set{D}$, so that $\rD{r}=\DR{R}$ for $r=\frac{e^R-1}{e^R+1}$ or $R=\ln\frac{1+r}{1-r}$.

If $\fol{}=\plfol{}$ is a \textsingfol{} and $x\in\manis{M}{E}$, we denote by $\leafu{x}$ the leaf of $\fol{}$ through $x$. Moreover, if $\leafu{x}$ is hyperbolic, we denote by $\phi_{x}\colon\set{D}\to\leafu{x}$ a uniformization of $\leafu{x}$ such that $\phi_x(0)=x$.

Given a Hermitian metric $\metm{M}$ on a complex manifold $\mani{M}$, we denote by $\dherm{M}{\cdot}{\cdot}$ the distance induced by $\metm{M}$. Similarly, we denote by $\dPC{\cdot}{\cdot}$ the Poincar\'{e} distance if one is given a Poincar\'{e} metric denoted $\metPC{}$.

Recall that $d^{c}=\frac{i}{2\pi}\left(\bhd-\hd\right)$ so that $\ddc=\frac{i}{\pi}\hd\bhd$.

\subsection*{Acknowledgments} The  author is supported by the Labex CEMPI (ANR-11-LABX-0007-01) and by the project QuaSiDy (ANR-21-CE40-0016).

  \section{Leafwise Poincar\'{e} metric\label{secPC}}

  Let $\fol{}=\plfol{}$ be a \textsingfol{}. Fix $\metm{M}$ a Hermitian metric on $\mani{M}$. Define
  \begin{equation}\label{defeta}\eta(x)=\sup\left\{\norm{\alpha'(0)}_{\metm{M}}\setst\alpha\colon\set{D}\to\leafu{x}~\text{holomorphic such that}~\alpha(0)=x\right\},\end{equation}
  where $\norm{v}_{\metm{M}}$ is the norm of $v\in T_x\leafu{x}$ with respect to the Hermitian metric $\metm{M}$. That is, $\norm{v}_{\metm{M}}=\sqrt{\metmp{M}{x}(v,v)}$.

  \begin{defn} A leaf $L$ of $\fol{}$ is called \emph{hyperbolic} if it is uniformized by the Poincar\'{e} disk $\set{D}$. The foliation $\fol{}$ is called \emph{hyperbolic} if all the leaves of $\fol{}$ are hyperbolic.
  \end{defn}
  
  We will be mostly interested in the case of hyperbolic foliations. Using Schwarz' Lemma, it is quite clear that $\leafu{x}$ is hyperbolic iff $\eta(x)<\infty$. In that case, $\eta(x)=\norm{\phi_x'(0)}_{\metm{M}}$, where $\phi_x\colon\set{D}\to\leafu{x}$ is a uniformization of $\leafu{x}$ such that $\phi_x(0)=x$. By pushing forward the Poincar\'{e} metric on $\set{D}$ by $\phi_x$, we obtain the so-called \emph{Poincar\'{e} metric $\metPC{}$ on the leaf $\leafu{x}$}. Since $\phi_x$ can be chosen up to pre-composition by a rotation, and since the biholomorphisms of $\set{D}$ are isometries of $\left(\set{D},\metPC{}\right)$, the Poincar\'{e} metric is intrinsically defined. It is quite clear that we have
  \begin{equation}\label{immpropeta}\eta^2\metPC{}=4\metm{M}.\end{equation}

    We recall from~\cite[Definition~3.1]{DNSII} the following.

\begin{defn} \label{defBrody}

We say that $\fol{}$ is \emph{Brody hyperbolic} if there exists a positive constant $c_0$ such that $\eta<c_0$ on $\manis{M}{E}$.

\end{defn}

Let us define the type of singularities we will study.

  \begin{defn}
    Let $\fol{}=\plfol{}$ be a \textsingfol{}. Near a singularity $a\in E$, there exists a vector field $X$ defining $\fol{}$. In coordinates $(z_1,\dots{},z_n)$ centered at $a$, we can write
    \[X(z)=\sum\limits_{j=1}^nF_j(z)\der{}{z_j}.\]
    The functions $F_j$ can be developed as a power series $F_j=\sum_{\alpha\in\set{N}^n}c_{\alpha,j}z^{\alpha}$. The \emph{1-jet} of $X$ at $a$ is defined in the chart $(U,z)$ as $X_1=\sum_{j=1}^n\sum_{\Cmod{\alpha}\leq1}c_{\alpha,j}z^{\alpha}\der{}{z_j}$. See for example~\cite[Chapter I]{Ilya} for more details. If the 1-jet of $X$ has an isolated singularity at $a$, we say that $a$ is a \emph{non-degenerate singularity of $\fol{}$}.
  \end{defn}

The following result is a direct consequence of a local study of non-degenerate vector fields in~\cite{MLN}. We denote by $\log^{\star}=1+\Cmod{\log}$ a log-type function and by $\dsing{M}{E}{x}$ the distance from $x\in\mani{M}$ to the singular set $\mani{E}$ with respect to the Hermitian metric $\metm{M}$.
  
  \begin{prop}[Dinh-Nguy\^{e}n-Sibony] \label{nondegMLN} Let $\fol{}=\plfol{}$ be a Brody hyperbolic compact \textsingfol{}. Suppose that all the singularities of $\fol{}$ are non-degenerate. Then, there exists a constant $c>1$ such that the Poincar\'{e} metric of $\fol{}$ satisfies
    \[c^{-1}\zlogz{\dsing{M}{E}{x}}\leq\eta(x)\leq c\zlogz{\dsing{M}{E}{x}},\qquad x\in\manis{M}{E}.\]
  \end{prop}

  See~\cite[Proposition~4.1]{Bac} for a proof, which is inspired by~\cite[Proposition~3.3]{DNSII}.

  \section{Leafwise heat diffusion operators\label{secDt}}

  We will now define one of the two heat diffusion operators we are interested in. Let $\fol{}=\plfol{}$ be a \texthypsingfol{}. The leafwise Poincar\'{e} metric $\metPC{}$ gives rise to its associated Laplacian $\LapP{}$ on leaves. More precisely, for $f\in\class{2}\left(\set{D}\right)$ or for $f\in\class{2}\left(\leafu{x}\right)$, $\LapP$ is defined by the following formula.
  \[\left(\LapP{}f\right)\metPC{}=\pi\ddc{}f=i\hd\bhd f\qquad\text{on}~\set{D}~\text{or}~\leafu{x}.\]

  For any fixed $x\in\manis{M}{E}$, we study the \emph{heat equation} on $\leafu{x}$
\begin{equation}\label{defheatker} \der{p(x,y,t)}{t}=\LapPv{y}p(x,y,t),\quad\underset{t\to0^+}{\lim}\,p(x,y,t)=\delta_x(y),\quad y\in\leafu{x},~t\in\set{R}_+.\end{equation}
Here, we have denoted by $\delta_x$ the Dirac mass in $x$ and by $\LapPv{y}$ the Laplacian $\LapP$ with respect to the variable $y$. The limit is taken in the sense of distribution. 

The \emph{heat kernel on $\leafu{x}$}, denoted by $p(x,y,t)$, is defined as the smallest positive solution of the heat equation~\eqref{defheatker} (see Chavel~\cite{Chavel}). By diffusing $p(x,y,t)$, we get a one-parameter family of operators $\left\{D_t\setst t\in\set{R}_+\right\}$ defined by
\begin{equation}\label{defDt}D_tf(x)=\int_{\leafu{x}}p(x,y,t)f(y)\metPC{}(y),\qquad x\in\manis{M}{E},\quad f\in L^{\infty}\left(\manis{M}{E}\right).\end{equation}
This family is a semi-group of positive contractions of $L^{\infty}\left(\manis{M}{E}\right)$ (see~\eqref{Dtsemigroup}). Its elements are called the \emph{leafwise heat diffusion operators}.
\begin{equation}\label{Dtsemigroup} D_0=\id;\quad D_t\mathbf{1}=\mathbf{1}\quad\text{and}\quad D_{t+s}=D_t\circ D_s,~\text{for}~t,s\in\set{R}_+,\end{equation}
where $\mathbf{1}$ is the function identically equal to $1$.

The same results hold on $\set{D}$ itself. We denote by $p_{\set{D}}(\theta,\zeta,t)$ the heat kernel on $\set{D}$. The heat diffusions on a leaf $\leafu{x}$ are then related to the one on $\set{D}$ by the formula
\begin{equation}\label{Dtunif} D_t\left(f\circ\phi_x\right)=\left(D_tf\right)\circ\phi_x,\qquad\text{on}~\set{D},~\text{for}~t\in\set{R}_+,~f\in L^{\infty}\left(\leafu{x}\right).\end{equation}

See~\cite[Proposition~2.7]{NguOse} for a proof.

  \section{Directed currents and harmonic measures\label{seccur}}

  We suppose that the reader is already familiar with the notion of currents, positive currents, harmonic currents. For a more detailed exposition, see~\cite{surVANG21}. For a new example, see~\cite{AlkReb}. Let $\fol{}=\plfol{}$ be a \texthypsingfol{}. A directed harmonic current on $\fol{}$ can be decomposed leafwise and transversally using the following lemma. Its proof was done in~\cite[Proposition~2.3]{DNS12} (see also~\cite[Proposition~2.5]{surVANG21}). We will need it to define the abstract heat diffusions.
  
  \begin{lem}\label{decharmcur} Let $T$ be a directed harmonic current on $\fol{}$. Consider $\set{U}\simeq\set{D}\times\set{T}$ a flow box the coordinates of which can be extended to a neighbourhood of its closure in $\mani{M}$. Then, there is a positive Radon measure $\nu$ on $\set{T}$ and for $\nu$-almost every $t\in\set{T}$ a harmonic function $h_t$ on $\set{D}$ such that
    \begin{itemize}
    \item The mass $\int_{\set{T}}\norm{h_t}_{L^1\left(\set{D}\right)}d\nu(t)$ is finite;
    \item For $\alpha\in\testformfol{1,1}$ compactly supported in $\set{U}$,
        \[T(\alpha)=\int_{\set{T}}\left(\int_{\set{D}}h_t(y)\alpha(y,t)\right)d\nu(t).\]
    \end{itemize}
      If moreover $T$ is positive, then for $\nu$-almost every $t\in\set{T}$, $h_t$ is positive on $\set{D}$.
  \end{lem}

  We have a notion of harmonicity for measures, characterizing their behaviour with respect to the Laplacians on the leaves.

\begin{defn} Let $\LapP$ be the aggregate of the leafwise Laplacians on the leaves of $\fol{}$. A finite positive Borel measure $\mu$ on $\manis{M}{E}$ is called \emph{harmonic} if
  \[\int_{\mani{M}}\LapP fd\mu=0,\qquad f\in\testfunfol{},\]
  where $\testfunfol{}=\testformfol{0,0}$ denotes the space of test functions.
\end{defn}

Harmonic measures and directed harmonic currents are strongly linked. More precisely, if $T$ is a directed harmonic current and if the measure
\begin{equation}\label{defmuT}\mu=T\wedge\metPC{}\quad\text{on}\quad\manis{M}{E},\end{equation}
is finite, then $\mu$ is harmonic. Conversely, if $\mu$ is a harmonic measure, then there exists a directed harmonic current $T$ such that~\eqref{defmuT}. See~\cite{DNS12} for an implicit proof. In any case, we call the total mass of $\mu$ defined by~\eqref{defmuT} the \emph{Poincar\'{e} mass} of $T$. It may be infinite near singularities. We will need the following results. The first one is a weak version of a theorem by Skoda in~\cite{Skoda}. The second is due to Dinh, Nguy\^{e}n and Sibony in~\cite{DNS12}.

\begin{prop}\label{propLelong} Let $T$ be a positive $\ddc$-closed $(0,0)$-current on a ball $r_0\set{B}$ and $\beta=\ddc\norm{z}^2$ be the standard K\"{a}hler form. Then, the quantity $r^{-2}\norm{T\wedge\beta}_{r\set{B}}$ is bounded for $r\in\intoc{0}{r_1}$, $r_1<r_0$.
\end{prop}
Here, $\norm{T\wedge\beta}_A$ denotes the mass of $T\wedge\beta$ as a positive measure on a measurable set $A$.

\begin{prop}\label{propextcurisolsing} Let $\fol{}=\plfol{}$ be a \textsingfol{} with isolated singularities and $T$ be a directed positive harmonic current on $\manis{M}{E}$. Then, $T$ has locally finite mass near any singularity $a\in E$ and can be extended through $a$ into a positive $\ddc$-closed current.
\end{prop}

The following criterion for the Poincar\'{e} mass of a $\ddc$-closed current to be finite will be very important to prove Theroem~\ref{impcor}.

\begin{prop} \label{nondegfinPCmass} Let $\fol{}=\plfol{}$ be a \textsingfol{} such that $\mani{E}$ is composed of isolated singularities. Suppose that in a neighbourhood of $a\in\mani{E}$, the Poincar\'{e} metric of $\fol{}$ satisfies $\eta(x)\geq c\zlogz{\dherm{M}{x}{a}}$ for a constant $c>0$. Then, any directed positive harmonic current on $\fol{}$ has locally finite Poincar\'{e} mass near $a$.
\end{prop}

\begin{pr} The proof was implicitly done in~\cite[Proposition 4.2]{DNS12}. Basically, we work in coordinates and with the standard K\"{a}hler form $\beta=\ddc\norm{z}^2$ of $\set{C}^k$. Propositions~\ref{propextcurisolsing} and~\ref{propLelong} ensure that $\norm{T\wedge\beta}_{r\set{B}}\leq Cr^2$. Then, we use~\eqref{immpropeta} and the proof ends with the following integration by parts.
  
    \[\begin{aligned}\int_{r\set{B}}T\wedge\metPC&=4\int_{r\set{B}}\frac{1}{\eta^2(z)}T\wedge\beta\leq\frac{4}{c^2}\int_{r\set{B}}\frac{1}{\norm{z}^2\ln^2\norm{z}}T\wedge\beta\\
  \int_{r\set{B}}T\wedge\metPC&\leq\frac{4}{c^2}\left[\frac{\norm{T\wedge\beta}_{r\set{B}}}{r^2\ln^2(r)}+\int_0^r\norm{T\wedge\beta}_{\rho\set{B}}\left(\frac{2}{\rho^3\ln^2\rho}+\frac{2}{\rho^3\ln^3\rho}\right)d\rho\right]\\
      &\leq\frac{8C}{c^2\Cmod{\ln r}}<+\infty.\end{aligned}\]
  \qed

\end{pr}

\section{Abstract heat diffusion\label{secSt}}

Let $\fol{}=\plfol{}$ be a \texthypsingfol{} and $\mu$ be a harmonic measure on $\fol{}$. We want to have a solution for an abstract heat equation. It will be given by Hille-Yosida's theorem. First recall some notions and facts of functional analysis. The reader can find proofs and exposition in Brezis~\cite{Brezis}.

Let $L$ be a Hilbert space and $\left\langle\cdot,\cdot\right\rangle$ be its scalar product. A linear operator $A\colon\Dom(A)\to L$, for $\Dom(A)\subset L$, is called \emph{monotone} if $\left\langle Au,u\right\rangle\geq0$, $u\in\Dom(A)$. It is called \emph{maximal monotone} if moreover for $f\in L$, there exists $u\in\Dom(A)$ such that $Au+u=f$. If $A$ is maximal monotone, then $\Dom(A)$ is dense in $L$ and the graph of $A$ is closed. A \emph{semi-group of contractions} is a family $S(t)\colon L\to L$, $t\in\set{R}_+$ satisfying
\begin{itemize}
  \item $S(t+t')=S(t)\circ S(t')$, $t,t'\in\set{R}_+$,
  \item $\norm{S(t)}\leq1$, $t\in\set{R}_+$.
\end{itemize}

  A maximal monotone operator $A$ gives rise to a semi-group of contractions having $A$ as infinitesimal generator by the following result.
  
\begin{thm}[Hille-Yosida] Let $A$ be a maximal monotone operator on a Hilbert space $L$. The equation
  \begin{equation}\label{eqSt}\der{u(t,\cdot)}{t}+Au(t,\cdot)=0,\quad\text{and}\quad u(0,\cdot)=u_0,\end{equation}
  for
  \begin{equation}\label{regSt}u\in\class{1}\left(\set{R}_+,L\right)\cap\class{0}\left(\set{R}_+,\Dom(A)\right)\end{equation}
  and $u_0\in\Dom(A)$, has a unique solution $u$. Moreover, $u$ is given by a semi-group of contractions $S(t)\colon L\to L$, $t\in\set{R}_+$ acting on $u_0$, \emph{i.e.} $u(t,\cdot)=S(t)u_0$.
\end{thm}

So, one will be able to solve the heat equation if one can show that the operator $-\LapP$ is maximal monotone. This work was done in~\cite[Proposition~5.6]{DNS12} (see also~\cite[Proposition~5.8]{surVANG21}) and has led to the following result.

\begin{prop}\label{existsSt} Let $\fol{}=\plfol{}$ be a \textsingfol{} endowed with its Poincar\'{e} metric $\metPC{}$. If $\mu$ is a harmonic measure on $\fol{}$, then $-\LapP$ is maximal monotone on $L^2(\mu)$.

  In particular, there exists a semi-group of contractions $S(t)$, $t\in\set{R}_+$ such that for $u_0\in\Dom\left(-\LapP\right)$, $u(t,\cdot)=S(t)u_0$ satisfies~\eqref{regSt} and~\eqref{eqSt}. The $S(t)$ are called the \emph{abstract heat diffusion operators}.
\end{prop}

The proof of~\cite{DNS12} relies on the study of bilinear forms defined through the operator $\LapP$ and another Laplace type operator $\tLapP$. We recall the definition of these bilinear forms and some of their properties we will need. Let $T$ be a directed harmonic current associated to $\mu$ by~\eqref{defmuT}. By Lemma~\ref{decharmcur}, we can write in a given flow box $\set{U}\simeq\set{D}\times\set{T}$
\[T(\alpha)=\int_{\set{T}}\left(\int_{\set{D}}h_t(y)\alpha(y,t)\right)d\nu(t),\]
where the $h_t$, for $t\in\set{T}$, are positive harmonic functions on $\set{D}$ and $\nu$ is a Radon measure. Define
\[\tLapP u=\LapP u+\left\langle h_t^{-1}\gr h_t,\gr u\right\rangle_{\metPC{}}.\]

Here, the gradient is defined with respect to the Poincar\'{e} metric so that $\left\langle\gr u,\xi\right\rangle_{\metPC{}}=du(\xi)$, for $u\in\testfunfol{}$ and $\xi$ tangent to the leaf. Note that $\Cmod{\gr u}=\Cmod{du}_{\PC{}}$. The gradient can be extended to an operator of domain $H^1(\mu)$ onto $L^2(\mu)$. The space $H^1(\mu)$ is the completion of $\testfunfol{}$ for the norm
\[\norm{u}_{H^1(\mu)}=\left(\norm{u}_{L^2(\mu)}^2+\norm{\gr u}_{L^2(\mu)}^2\right)^{\frac{1}{2}}.\]
We have $\Dom\left(-\LapP\right)\subset H^1(\mu)$. \emph{A priori}, the definition of $\tLapP$ depends on the choice of flow boxes, but the uniqueness of the functions $h_t$ and the current $T$ ensure that it is a global operator. Define the bilinear forms
\[q(u,v)=-\int(\LapP u)v\,d\mu,\qquad \tilde{q}(u,v)=-\int(\tLapP u)v\,d\mu,\]
for $u,v\in\testfunfol{}$. We summarize the properties of $q$ and $\tilde{q}$ we will need in Section~\ref{secidheats} in the following statement. The reader can find proofs in~\cite[Lemmas~5.4 and~5.5]{DNS12} (see also~\cite[Lemmas~5.4 and~5.5, Remark~5.7]{surVANG21}).
\begin{prop}\label{factsq} Let $\mu$ be a harmonic measure. The bilinear forms $q$ and $\tilde{q}$ can be extended continuously to $H^1(\mu)\times H^1(\mu)$. Moreover,
  \begin{enumerate}
  \item \label{q=tq} For $u\in H^1(\mu)$, $q(u,u)=\tilde{q}(u,u)$.
  \item \label{formtq} For $u\in\Dom\left(-\LapP\right)$ and $v\in H^1(\mu)$, we have
    \[\tilde{q}(u,v)=\int\left\langle\gr u,\gr v\right\rangle_{\metPC{}}d\mu=\int i\hd u\wedge\bhd v\wedge T.\]
  \end{enumerate}
\end{prop}

\section{Coincidence of the heat diffusions : proof of Theorem~\ref{impcor}\label{secidheats}}

The aim of this section is to compare the diffusion operators $S(t)$ defined in Section~\ref{secSt} and $D_t$ defined in Section~\ref{secDt}, for any $t\in\set{R}_+$. We begin by stating a very abstract criterion for both heat diffusions to be equal. 

\begin{thm}\label{mainthm} Let $\fol=\plfol{}$ be a \textsingfol{} with isolated singularities. Consider on $\mani{M}$ the leafwise Poincar\'{e} metric $\metPC{}$ and a Hermitian metric $\metm{M}$. Denote by $\eta$ the ratio function defined by~\eqref{defeta}. Let $\mu$ be a harmonic measure on $\fol{}$. Assume that
  \begin{enumerate}[label=(H\arabic*)]
  \item\label{hypmtbndeta} $\eta$ is locally bounded from above on $\manis{M}{E}$ and in the neighbourhood of any $a\in\mani{E}$, $\eta(x)\leq\upb_a\left(\dherm{M}{x}{a}\right)$ for some function $\upb_a\colon\set{R}_+^*\to\set{R}_+^*$ such that $\upb_a$ is continuous and $\frac{1}{\upb_a}$ is not integrable near $0$;
  \item\label{hypmtmember} Any measurable function $u$ on $\manis{M}{E}$ satisfying
    \begin{enumerate}
    \item $\norm{u}_{L^{\infty}}<\infty$,
    \item $\norm{\Cmod{du}_{\PC}}_{L^{\infty}}<\infty$,
    \item $\norm{\LapP u}_{L^{\infty}}<\infty$,
    \item for every $a\in\mani{E}$, $\underset{\eps\to0}{\lim}\,\norm{\restriction{u}{\left\{\dherm{M}{x}{a}<\eps\right\}}}_{L^{\infty}}=0$,
    \end{enumerate}
    belongs to $H^1(\mu)$.
  \end{enumerate}
  
  Then, the abstract heat diffusions operators and the leafwise heat diffusions operators coincide. 
\end{thm}

\begin{rem}\label{analogabstractthmVANG} This theorem is analogous to a result by Nguy\^{e}n~\cite[Theorem~5.17]{surVANG21}. In comparison, our hypothesis~\ref{hypmtbndeta} is stronger than his, and our hypothesis~\ref{hypmtmember} is weaker than his. Actually, he uses the first one in a context where he has far better. Therefore, he is forced to suppose quite specific singularities and currents in order to satisfy the second one. Namely, he has to work with hyperbolic singularities and a current that does not give mass to any leaf to have a vanishing Lelong number. Here, we have relaxed the second hypothesis to apply the theorem to more singular measures, and strengthened the first one in consequence. In particular, our assumption~\ref{hypmtbndeta} allows us to obtain the result of Lemma~\ref{factepsinf} below.
\end{rem}

Before proving Theorem~\ref{mainthm}, we will need two lemmas. The first one is new and will be a crucial ingredient of our proof. The second one can be found in~\cite{surVANG21}.

\begin{lem}\label{factepsinf} We keep the notations of Theorem~\ref{mainthm} and suppose hypothesis~\ref{hypmtbndeta}. Let $u_0\in\testfunfol{}$ and $u(t,\cdot)=D_tu_0$. Then, for $t>0$ and $a\in E$, $\underset{\eps\to0}{\lim}\norm{\restriction{u}{\{\dherm{M}{x}{a}<\eps\}}(t,\cdot)}_{L^{\infty}}=0$.
\end{lem}

\begin{pr} Let $a\in E$ and $\delta$ be such that there exists a Hermitian chart $\set{U}$ around $a$ of radius $\delta$. Without loss of generality, we can assume that $\metm{M}$ is the standard K\"{a}hler metric on $\set{U}\simeq\delta\set{B}$. Shrinking the chart if necessary, we also suppose that $\delta<\dherm{M}{\supp u_0}{a}$. Fix $t>0$, take $\eps<\delta$ and $x\in\eps\set{B}$. For $y\in\leafu{x}\cap\supp u_0$, we have by~\eqref{immpropeta},
  \[\dPC{x}{y}\geq2\int_{\eps}^{\delta}\frac{dt}{F_a(t)}=R_{\eps}\underset{\eps\to0}{\to}+\infty,\]
  since $F_a$ is not integrable near $0$. It follows that $\phi_x^{-1}\left(\supp u_0\right)\subset\wo{\set{D}}{\DR{R_{\eps}}}$. In other words, in the uniformization of $\leafu{x}$ centered at $x$, the support of $u_0$ is sent uniformly to the boundary when $x$ goes to the singularity. By~\eqref{Dtunif} and~\eqref{defDt}, we get
  \[\begin{aligned} D_tu_0(x)=u(t,x)&=\int_{\set{D}}p_{\set{D}}(0,y,t)u_0\left(\phi_x(y)\right)\metPC{}(y),\\
      \Cmod{u(t,x)}&\leq\norm{u_0}_{L^{\infty}}\int_{\set{D}}\indic{\wo{\set{D}}{\DR{R_{\eps}}}}(y)p_{\set{D}}(0,y,t)\metPC{}(y),\end{aligned}\]
  where $\indic{\wo{\set{D}}{\DR{R_{\eps}}}}$ denotes the characteristic function of $\wo{\set{D}}{\DR{R_{\eps}}}$.  This bound is independent on $x$ such that $\dherm{M}{x}{a}<\eps$. By Lebesgue's dominated convergence theorem and the fact that $R_{\eps}\underset{\eps\to0}{\to}+\infty$, the right hand side tends to $0$ when $\eps$ goes to $0$. This implies the desired conclusion.
  \qed
\end{pr}

\begin{lem}\label{factslemma} We keep the notations and hypothesis of Lemma~\ref{factepsinf}. Suppose moreover that $u_0$ is positive. Then $u$ is measurable and satisfies
    \begin{enumerate}
  \item\label{factmembership} \begin{enumerate}
    \item\label{factinf} $\norm{u(\cdot,\cdot)}_{L^{\infty}}<\infty$,
    \item\label{factdinf} For $t>0$, $\norm{\Cmod{du(t,\cdot)}_{\PC}}_{L^{\infty}}<\infty$,
    \item\label{factLapinf} For $t>0$, $\norm{\LapP{u}(t,\cdot)}_{L^{\infty}}<\infty$;
    \end{enumerate}
  \item\label{factheateq} $u\in\class{1}\left(\set{R}_+^*,L^2(\mu)\right)$ and for $t\in\set{R}_+^*$, $\der{u(t,\cdot)}{t}-\LapP u(t,\cdot)=0$;
  \item\label{factlim} $\underset{t\to0}{\lim}\,u(t,\cdot)=u_0$ in $L^2(\mu)$.
  \end{enumerate}
\end{lem}

See~\cite[Facts (5.15), (5.16) and (5.17)]{surVANG21} for a proof relying on a precise estimate of the heat kernel.

\emph{Proof of Theorem~\ref{mainthm}.} We follow the proof of~\cite[Theorem~5.17]{surVANG21}. Let $u_0\in \Dom\left(-\LapP\right)$ and $t\in\set{R}_+$. We have to show that
\begin{equation}\label{eqheateq}D_tu_0=S(t)u_0.\end{equation}

Since the operators $D_t$ and $S(t)$ are positive contractions and $\testfunfol{}$ is dense in $\Dom\left(-\LapP\right)$, it is enough to show~\eqref{eqheateq} for a non-negative $u_0\in\testfunfol{}$. Take such a $u_0$. Denote by $u(t,\cdot)=D_tu_0$ and $U(t,\cdot)=S(t)u_0$. By Lemmas~\ref{factepsinf} and~\ref{factslemma} and hypothesis~\ref{hypmtmember}, $u(t,\cdot)$, $U(t,\cdot)$ and $\tilde{u}(t,\cdot)=U(t,\cdot)-u(t,\cdot)$ belong to $H^1(\mu)$ for a fixed $t$. Morevoer, they are $\class{1}\left(\set{R}_+^*,L^2(\mu)\right)$. The heat equations~\eqref{eqSt} and~\eqref{factheateq} in Lemma~\ref{factslemma} satisfied by $S(t)$ and $D_t$ imply that
\[\frac{1}{2}\der{}{t}\left(\norm{\tilde{u}(t,\cdot)}_{L^2(\mu)}^2\right)=\left\langle\LapP\tilde{u}(t,\cdot),\tilde{u}(t,\cdot)\right\rangle=-q\left(\tilde{u}(t,\cdot),\tilde{u}(t,\cdot)\right).\]

By an approximation of $\tilde{u}(t,\cdot)$ by elements of $\testfunfol{}$ and Proposition~\ref{factsq}, we get $q\left(\tilde{u}(t,\cdot),\tilde{u}(t,\cdot)\right)=\tilde{q}\left(\tilde{u}(t,\cdot),\tilde{u}(t,\cdot)\right)\geq0$ (see~\cite[p.~50-51]{surVANG21}). Hence, $\norm{\tilde{u}(t,\cdot)}_{L^2(\mu)}^2$ is decreasing and goes to $0$ for $t\to0$ by the boundary conditions. This implies that it is identically zero and $U(t,\cdot)=u(t,\cdot)$.
\qed

\emph{Proof of Theorem~\ref{impcor}.} We are inspired by the proof of~\cite[Proposition~5.22]{surVANG21}. Let $T$ be a directed positive harmonic current on $\fol{}$. By Proposition~\ref{nondegMLN}, we have
\begin{equation}\label{bndMLN}c^{-1}\zlogz{\dsing{M}{E}{x}}\leq\eta(x)\leq c\zlogz{\dsing{M}{E}{x}},\quad x\in\manis{M}{E}.\end{equation}
Then, Proposition~\ref{nondegfinPCmass} ensures that the measure $\mu$ associated to $T$ by~\eqref{defmuT} is finite. Hence, it is a harmonic measure. We want to apply Theorem~\ref{mainthm}. Hypothesis~\ref{hypmtbndeta} is clear by~\eqref{bndMLN}. It remains to check hypothesis~\ref{hypmtmember}. Let $u$ be a measurable function on $\manis{M}{E}$ such that
\begin{equation} \label{hypmember} \norm{u}_{L^{\infty}}<\infty,\quad\norm{\Cmod{du}_{\PC{}}}_{L^{\infty}}<\infty,\quad\norm{\LapP u}_{L^{\infty}}<\infty,\quad\underset{\eps\to0}{\lim}\norm{\restriction{u}{\{\dherm{M}{x}{a}<\eps\}}}_{L^{\infty}}=0,~a\in E.\end{equation}
We need to show that $u$ belongs to $H^1(\mu)$. Using the compactness of $\mani{M}$, we can find a finite open covering $\mathcal{U}=\left(U_p\right)_{p\in I}$ by
  \begin{enumerate}
  \item \label{pnotinE} flow boxes $U_p\simeq\set{D}\times\set{T}$ such that $\dherm{M}{U_p}{E}>0$,
  \item \label{pinE} Hermitian charts $U_a\simeq\set{B}$ for $a\in E$ with $U_a\cap E=\{a\}$.
  \end{enumerate}

  Since $\mu$ is finite, using a partition of the unity, we can suppose that $u$ is compactly supported in a single $U_p$ or $U_a$.

  \emph{Case~\eqref{pnotinE}.} We need to find $u_{\eps}\in\testfunfol{}$ such that $\norm{u-u_{\eps}}_{H^1(\mu)}<\eps$. We will follow the three following steps.
  \begin{enumerate}[label=(\roman*)]
  \item Let $\chi^{(1)}\colon\set{C}\to\intcc{0}{1}$ be a smooth function supported in $\adh{\set{D}}$ such that $\int_{\set{C}}\chi^{(1)}d\Leb=1$. Define $\chi^{(1)}_{\delta}(z)=\delta^{-2}\chi(\delta z)$ for $\delta\in\intoo{0}{1}$. Let $u_1$ be the leafwise convolution of $u$ with $\chi^{(1)}_{\delta}$.
  \item Let $v$ be a bounded continuous function such that $\norm{du_1-v}_{L^2(m)}<\eps$ and $w$ be a bounded continuous function such that $\norm{u_1(0,t)-w(t)}_{L^2(\nu')}<\eps$. Here, $m$ and $\nu'$ are measures respectively on $\set{D}\times\set{T}$ and $\set{T}$ that our computation will specify. Recall that we can find such $v$ and $w$ if $m$ and $\nu'$ are Radon measures by Luzin's Theorem. Define
    \[u_2(z,t)=w(t)+\int_0^1v(\tau z,t)\cdot zd\tau.\]
  \item Let $u_{\eps}=u_3$ be the leafwise convolution of $u_2$ with $\chi^{(1)}_{\delta}$.
  \end{enumerate}
  
  Note first that by~\eqref{bndMLN} and~\eqref{immpropeta}, there is a constant $C>1$ such that $C^{-1}\norm{dz}^2\leq\metPC{}\leq C\norm{dz}^2$ on $U_p$. It follows that the first and third steps will give $\norm{u_1-u}_{H^1(\mu)}<\frac{\eps}{3}$ and $\norm{u_3-u_2}_{H^1(\mu)}<\frac{\eps}{3}$, for sufficiently small $\delta$. Moreover, it is clear that $u_2$ is continuous. Thus, $u_3\in\testfunfol{}$. So, it remains to prove that $\norm{u_1-u_2}_{H^1(\mu)}<\frac{\eps}{3}$. Since $u_1$ is leafwise smooth, we have
  \[u_1(z,t)=u_1(0,t)+\int_0^1du_1(\tau z,t)\cdot zd\tau\]
  everywhere on $(z,t)\in\set{D}\times\set{T}$. Recall that $\mu$ is given by $T\wedge\metPC{}$ and $T$ by Lemma~\ref{decharmcur}. Decompose $u_1-u_2$ into the transversal to $0$ part and the integral part, and apply the Cauchy-Schwarz inequality to the second part. We get
  \[\norm{u_1-u_2}_{L^2(\mu)}\leq\norm{w-u_1(0,\cdot)}_{L^2(\nu')}+\norm{v-du_1}_{L^2(m_1)},\]
where $\nu'(t)=\norm{h_t}_{L^1(\metPC{})}\nu(t)$ and $m_1=\lambda_*\left(\Leb\times\mu\right)$ with $\lambda\colon(\tau,z,t)\mapsto(\tau z,t)$. Moreover, $du_2=v$ and since $C^{-1}\norm{dz}^2\leq\metPC{}\leq C\norm{dz}²$, we have
\[\norm{\gr u_1-\gr u_2}_{L^2(\mu)}\leq C\norm{v-du_1}_{L^2(\mu)}.\]

It is clear that $m=m_1+\mu$ is finite, hence is a Radon measure. So is $\nu'$ by Lemma~\ref{decharmcur}. Choose $v$ a bounded continuous function such that $\norm{du_1-v}_{L^2(m)}<\alpha\eps$ and $w$ a bounded continuous function such that $\norm{u_1(0,\cdot)-w}_{L^2(\nu')}<\alpha\eps$, for $\alpha\leq\frac{1}{3}\left(4+C^2\right)^{-1/2}$. Then $\norm{u_2-u_1}_{H^1(\mu)}<\frac{\eps}{3}$ and therefore, $\norm{u-u_{\eps}}_{H^1(\mu)}<\eps$. This concludes Case~\eqref{pnotinE}.

\emph{Case~\eqref{pinE}.} We work in coordinates $U_a\simeq\set{B}\subset\set{C}^k$ centered at $a$. Let $\chi\colon\set{C}^k\to\intcc{0}{1}$ be a smooth function such that $\chi(z)=1$ on $\frac{1}{2}\set{B}$ and $\chi(z)=0$ on $\wo{\set{C}^k}{\set{B}}$. For $\eps\in\intoo{0}{1}$, consider the function
\[v_{\eps}(z)=\left(1-\chi\left(\eps^{-1}z\right)\right)u(z),\qquad z\in\set{B}.\]

By construction, we have $v_{\eps}=0$ on $\frac{\eps}{2}\set{B}$ and $v_{\eps}$ satisfies~\eqref{hypmember}. By the same arguments as in Case~\eqref{pnotinE}, we get that $v_{\eps}\in H^1(\mu)$. So, it will be sufficient to prove that $\lim_{\eps\to0}\norm{u-v_{\eps}}_{H^1(\mu)}=0$. Moreover, $u-v_{\eps}$ is supported in $\eps\set{B}$ and $\Cmod{u-v_{\eps}}\leq\Cmod{u}$ everywhere. Since $\mu\left(\{a\}\right)=0$,
\[\int_{\mani{M}}\Cmod{u-v_{\eps}}^2d\mu\leq\norm{u}_{L^{\infty}}^2\int_{\eps\set{B}}d\mu\underset{\eps\to0}{\to}0.\]

Next, we adress the gradient-part of the $H^1(\mu)$-norm. We have
\[\int_{\mani{M}}\Cmod{\gr\left(u-v_{\eps}\right)}^2d\mu=\int_{\eps\set{B}}i\hd\left(u(z)\chi\left(\eps^{-1}z\right)\right)\wedge\bhd\left(u(z)\chi\left(\eps^{-1}z\right)\right)\wedge T(z).\]
Up to multiplying by constants, the right hand side is lower than the sum of these three terms
\[I_1=\int_{\eps\set{B}}i\hd u\wedge\bhd u\wedge T,\quad I_2=\eps^{-1}\int_{\wo{\eps\set{B}}{\frac{\eps}{2}\set{B}}}\Cmod{u}\Cmod{du}T\wedge\beta,\quad I_3=\eps^{-2}\int_{\wo{\eps\set{B}}{\frac{\eps}{2}\set{B}}}\Cmod{u}^2T\wedge\beta,\]
where $\beta=\ddc\norm{z}^2$ is the standard K\"{a}hler metric. Since $\mu\left(\{a\}\right)=0$, we obtain
\[I_1=\int_{\eps\set{B}}\Cmod{\gr u}^2d\mu\leq\norm{\Cmod{du}_{\PC{}}}^2\mu\left(\eps\set{B}\right)\underset{\eps\to0}{\to}0.\]

Recall that $\Cmod{du}=\frac{\Cmod{du}_{\PC}}{\eta}=\frac{\Cmod{\gr u}}{\eta}$. By~\eqref{bndMLN} and the equivalence of $\metm{M}$ and the standard norm, we get
\[I_2\leq \frac{C\norm{u}_{L^{\infty}}}{\log^{\star}(c\eps)}\norm{\Cmod{du}_{\PC}}_{L^{\infty}}\eps^{-2}\norm{T\wedge\beta}_{\eps\set{B}}.\]
Then, Proposition~\ref{propLelong} implies that $I_2\to0$ when $\eps\to0$. Finally,
\[I_3\leq\norm{\restriction{u}{\{\dherm{M}{x}{a}<\eps\}}}_{L^{\infty}}^2\eps^{-2}\norm{T\wedge\beta}_{\eps\set{B}}.\]
Thus, Proposition~\ref{propLelong} and~\eqref{hypmember} ensure that $I_3\to0$ when $\eps\to0$. It follows that $\norm{u-v_{\eps}}_{H^1(\mu)}$ goes to $0$ when $\eps\to0$. Hence, $u\in H^1(\mu)$ and Theorem~\ref{impcor} is a consequence of Theorem~\ref{mainthm}.
\qed

\begin{rem} In fact, $I_1$ is estimated exactly the same way as~\cite[p. 56-57]{surVANG21} and $I_2$ similarly. The new point here is the estimate on $I_3$ which is based on Lemma~\ref{factepsinf}.
\end{rem}

\end{document}